# VECTOR OPTIMIZATION BY TWO OBJECTIVE FUNCTIONS

Bosov[1], A. A., Kodola[1,2], G. N., Savchenko, L. N.

[1] Dniepropetrovsk National University of Railway Transport named after academician V.A. Lazarjan,
2 academician V.A. Lazarjan Street, Dnepropetrovsk, 49010, Ukraine
*e-mail: AABosov@i.ua*
[2] *e-mail:Gylina@i.ua*

Recently wide application in engineering-economic problems was received with problems of vector optimization. Development of methods of the decision of these problems it is executed in works A. Messac and others. Complexity of the offered methods consists in construction of an aggregate objective function (AOF). In the given work an algorithm of the solution of a vector optimization problem is suggested carry out by use analytical representation of Pareto cone.

**Introduction**

In his work [1] Euler L. pointed to the fact that any problem solved leads us to a task of maximum or minimum.

A mathematical formulation of the classical optimization problems is

$$f(x) \to \min$$

where

$$x \in X \subset E_n,$$

where $E_n$ – $n$-dimensional Euclidean space.

Such kind of a statement of the task is highly general.

Different functions $f(x)$ and feasible areas of $X$ allow us to consider a problem of existence of solutions and methods of solving the optimization problem [2].

However, there are too many engineering and economical problems that cannot be solved by classical methods.

The main feature of such kind of problems is rational use of resources.

The rationality of their use is evaluated by objective functions (factors) $f_i(x)$, $i = \overline{1,k}$, and each one must be minimized when $X$ is defined.

A formal note of this kind of problems is

$$\begin{pmatrix} f_1(x) \\ f_2(x) \\ ... \\ f_n(x) \end{pmatrix} \to \min$$

where $x \in X$.

The main advantage of such statement of problem is the fact that it makes the problem demonstrable and we can formulate the solution rule of the problem.

However, the solution of such kind is ambiguous.

In mathematical literature such problems are formulated with a use of the binary relations [3,4].

# VECTOR OPTIMIZATION BY TWO OBJECTIVE FUNCTIONS

The Pareto`s binary relation is used as a selection criterion of solution variants for formulated problem.

There are many literature devoted to such kind of problems. Let us consider the works [5-6], where we can find a detailed survey of literature devoted to the vector optimization problem and the methods of its solving.

As a rule, vector optimization methods lead to construction of an aggregate objective function (AOF) [7]. The construction method of AOF in a form of linear (weighting) combination of the objective functions is the most wide-spread one. And the main weak point of the method is determination of the weighting coefficients.

The necessary conditions for aggregate objective function that allow to determine the full Pareto`s surface with using one of the next methods were received in works [7-9].

The NBI (Normally Boundary Intersection) method described in work [10] is also widely spread. However, when we use the NBI method we can get the non-optimum points by Pareto and also local points that need to be filtered. A new NC method (Normal Constraint), which was got in works [11-12] ensures the full representation of Pareto`s boundary, though it also can receive non-optimal points, but with not so great probability than the NBI method.

In this paper a new method of construction of Pareto`s boundary is suggested. The method is close to the one that you can find in works by Messac A.

First let us consider one of the mapping construction methods that is underlain construction method of Pareto`s solution set of vector optimization problem.

**Mapping construction**

Let $F(x)$ be defined, $x \in X \subseteq R_n$, as

$$F(x) = (F_1(x), F_2(x), ..., F_k(x)).$$

A set $Y \subseteq R_k$ is called mapping or image of the set $X$, if

$$Y = \{y \in R_k : y = F(x), x \in X\}. \tag{1}$$

Let vector $u \in R_n$ be of unite length and we consider a ray

$$x = u \cdot \tau. \tag{2}$$

Let interval $[\underline{\tau}(u), \overline{\tau}(u)]$ be such an interval that when $\tau \in [\underline{\tau}(u), \overline{\tau}(u)]$ ray points (2) belong to the set $X$. Then the set (1) can be determined as

$$Y = \{y \in R_k : y = F(u \cdot \tau); |u| = 1, \tau \in [\underline{\tau}(u), \overline{\tau}(u)]\}. \tag{3}$$

The vector $u$ can be determined in a form of

$$\begin{aligned}
u_1 &= \cos\varphi_1; \\
u_2 &= \sin\varphi_1 \cos\varphi_2; \\
&\cdots\cdots\cdots \\
u_i &= \sin\varphi_1 \sin\varphi_2 ... \sin\varphi_{i-1} \cos\varphi_i; \\
&\cdots\cdots\cdots \\
u_n &= \sin\varphi_1 \sin\varphi_2 ... \sin\varphi_{n-2} \sin\varphi_{n-1},
\end{aligned} \tag{4}$$

where $\varphi_1 \in [0, \pi]$; $\varphi_i \in [0, 2\pi]$, $i \geq 2$.

# VECTOR OPTIMIZATION BY TWO OBJECTIVE FUNCTIONS

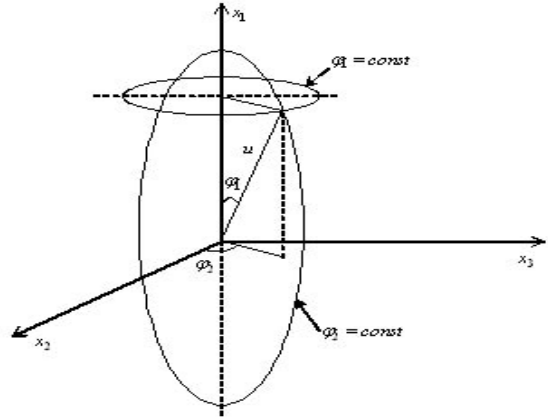

Figure 1.

You can see an interpretation of relations (4) for 3D-space in a fig.1. A line $\varphi_1 = const$ is known as latitude, and $\varphi_2 = const$ is a meridian of sphere with unit radius.

Note that relation (3) can be applied, when set $X$ is convex.

If $X$ is not convex, then we introduce its indicator

$$I_X(x) = \begin{cases} 1, & \text{if } x \in X; \\ 0, & \text{if } x \notin X. \end{cases}$$

In this case in relation (3) interval $\tau \in [\underline{\tau}(u), \overline{\tau}(u)]$ must be converted to a set

$$T(u) = \{\tau : I_X(u \cdot \tau) = 1\}$$

and then

$$Y = \{y \in R_k : y = F(u \cdot \tau); |u| = 1, \tau \in T(u)\}. \tag{5}$$

For numerical solution let the set $X$ be

$$X = \{x \in R_n : h_i(x) \leq 0; i = \overline{1,m}\}. \tag{6}$$

In this case we introduce function

$$H(x) = \max_{1 \leq i \leq m} \{h_i(x)\}, \tag{7}$$

then it obtains

$$H(x) = \begin{cases} \leq 0, & \text{if } x \in X; \\ > 0, & \text{if } x \notin X, \end{cases}$$

and the set $T(u)$ can be determined as

$$T(u) = \{\tau \in R_1 : H(u \cdot \tau) \leq 0\}. \tag{8}$$

Note, that the roots of the equation $H(u \cdot \tau) = 0$ let us determine the points of intersection of ray (2) with boundary of the set $X$, if the set is closed (see fig.2).

# VECTOR OPTIMIZATION BY TWO OBJECTIVE FUNCTIONS

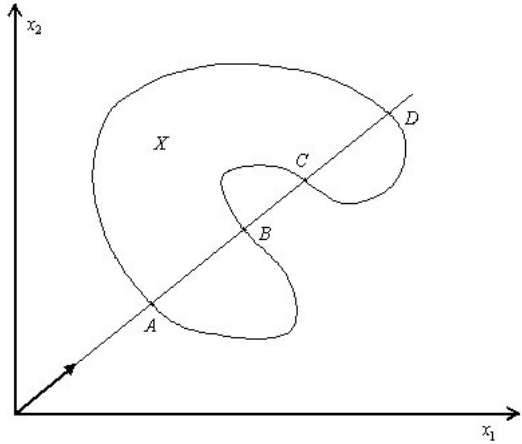

Figure 2. An intersection of the ray (2) with the set $X$.

If the boundary of the set $X$ does not belong to the set $X$, then the roots of the equation $H(u \cdot \tau)$ let us find the boundaries of the set $T(u)$ for exhaustive search $\tau$ during the construction of the intersection of the ray (2) with $X$.

Thus, when the vector $u$ is fixed, we can build next set

$$X(u) = \{x \in X : x = u \cdot \tau; H(u \cdot \tau) \leq 0\}.$$

From the example in fig. 2 we can see that this set consists of two intervals $[A, B]$ and $[C, D]$.

The set $X(u)$ is always one-dimensional set as an image of one-dimensional set $T(u)$.

Let the image of the set $X(u)$ be determined with $Y(u) \subseteq Y$ during the mapping by the vector-function $F(x)$.

The set $Y(u)$ is also one-dimensional set and in general case it represents itself some curve in a space $R_k$.

It is obvious, that if $F(x)$ is

$$F(x) = c \cdot x + b$$

where $C$ − matrix $k \times n$, then $Y(u)$ will be an interval or a set of intervals, because

$$F(\alpha x_1 + (1-\alpha)x_2) = C(\alpha x_1 + (1-\alpha)x_2) + b =$$
$$= \alpha(Cx_1 + b) + (1-\alpha)(Cx_2 + b) = \alpha F(x_1) + (1-\alpha)F(x_2)$$

when $\alpha \in [0,1]$, i. e. an interval $[x_1, x_2] \in R_n$ turns into an interval $[F(x_1), F(x_2)] \in R_k$.

Example 1. This example was taken from the work [13]. The set $X$ is

$$X = \{(x_1, x_2) \in R_2 : x_1 + x_2 \geq 1, x_1, x_2 \geq 0\}.$$

Using denotations (6) we get

$$h_1(x) = 1 - x_1 - x_2;$$
$$h_2(x) = -x_1;$$
$$h_3(x) = -x_2.$$

# VECTOR OPTIMIZATION BY TWO OBJECTIVE FUNCTIONS

$$F_1(x) = \frac{-x_1 + 0.5}{x_1 + x_2 - 0.75};$$

$$F_2(x) = \frac{-x_2 + 0.5}{x_1 + x_2 - 0.75}.$$

Let $u_1 = \cos\varphi$; $u_2 = \sin\varphi$, then

$$T(u) = T(\varphi) = \left[\frac{1}{\cos\varphi + \sin\varphi}, \infty\right],$$

And the set $Y$ is

$$Y = \left\{(y_1, y_2) \in R_2 : y_1 = \frac{-\cos\varphi \cdot \tau + 0.5}{(\cos\varphi + \sin\varphi)\tau - 0.75}; y_2 = \frac{-\sin\varphi \cdot \tau + 0.5}{(\cos\varphi + \sin\varphi)\tau - 0.75}; \tau \in T(\varphi), \varphi \in [0, \frac{\pi}{2}]\right\}.$$

The set $Y$ is shown in a fig.3.

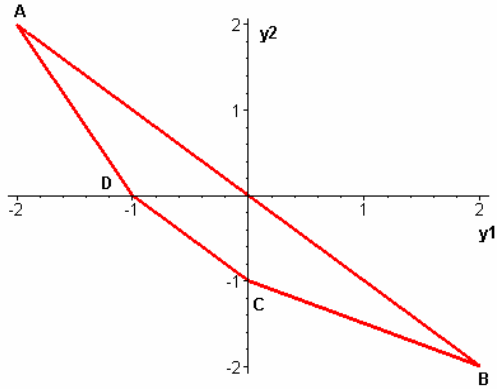

Figure 3. Geometrical mapping of the set $Y$ from the example 1.

Note that the interval $AB$ is an image of the interval $x_1 + x_2 = 1$ when $x_1, x_2 \geq 0$, and infinitely distant points of the set $X$ in the set $Y$ are represented with the interval $DC$.

**Method for construction Pareto frontier**

Let us consider vector optimization problem by two objective functions $f_1(x)$ and $f_2(x)$, $x \in R_n$ for convenience of geometrical interpretation. Each function must be minimized and a formal note of minimization is

$$\begin{pmatrix} f_1(x) \\ f_2(x) \end{pmatrix} \to \min,$$

where

$$x \in X \subseteq R_n.$$

Let

$$y_1(x) = f_1(x);$$
$$y_2(x) = f_2(x),$$

# VECTOR OPTIMIZATION BY TWO OBJECTIVE FUNCTIONS

then we can map the set $X$ into the set $Y \subseteq R_2$ and the original problem is

$$\begin{pmatrix} y_1 \\ y_2 \end{pmatrix} \to \min \qquad (9)$$

if $y \in Y$.

Remind, that the solution of the problem (9) is the set $Y_*$, whose points are incomparable due to Pareto.

Let $y_* \in Y_*$, and $K$ – cone, whose top is at the point $y_*$, when for every $y \in K$ satisfying next condition

$$\begin{cases} y_1 \leq u_1 \cdot t; \\ y_2 \leq u_2 \cdot t, \end{cases}$$

where $u_1$, $u_2$ – components of unite vector $u$, such that

$$\begin{cases} y_{1*} = u_1 \cdot t; \\ y_{2*} = u_2 \cdot t, \end{cases}$$

then

$$K \cap Y = \{y_*\}. \qquad (10)$$

Condition (10) is necessary and sufficient for $y_*$ to belong to the solution of the task (9).

Geometrical interpretation of the condition (10) is shown in a fig. 4.

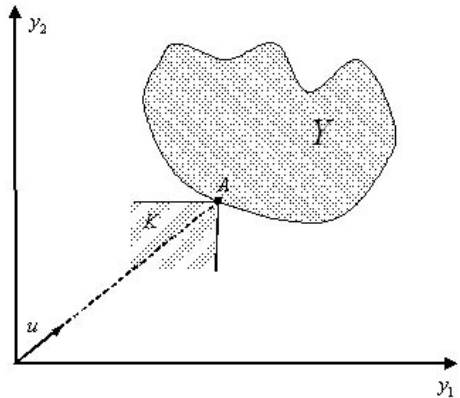

Figure 4. Geometrical interpretation of the condition (10).

Hereinafter we will suggest that the set $Y$ is obtained from the relation (6) by the next method:

$$Y = \{y \in R_2 : h_i(y) \leq 0, i = \overline{1, k}\} \qquad (11)$$

Let $y(1)$ be the solution of a task

$$y_1 \to \min,$$

if $y \in Y$, and $y(2)$ is the solution of this task under the same condition

$$y_2 \to \min.$$

# VECTOR OPTIMIZATION BY TWO OBJECTIVE FUNCTIONS

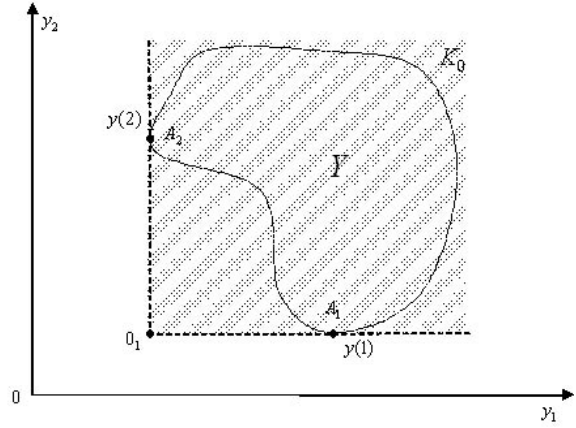

Figure 5. Geometrical interpretation of the tasks $y_1 \to \min$; $y_2 \to \min$.

Hereinafter we will transpose the origin of coordinates to the point $0_1$, and take $0_1 A_1$ and $0_1 A_2$ as the axes, and the area $Y \in K_0$ (see fig. 5). In this coordinate system "old" coordinates are represented by "new" ones $\tilde{y}_1$ and $\tilde{y}_2$

$$\begin{cases} y_1 = \tilde{y}_1 + y(1); \\ y_2 = \tilde{y}_2 + y(2). \end{cases} \quad (12)$$

Let us denote vector with components $(\tilde{y}_1, \tilde{y}_2)$ by $(y_1, y_2)$ for convenience that equals hypothesis, which we can understand from a fig. 6.

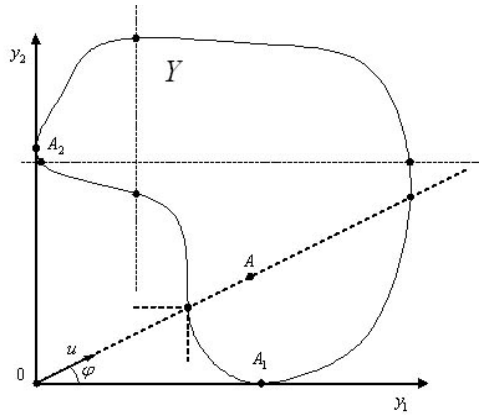

Figure 6. Geometrical representation of the area $Y$ after transformation (12)

Using (4) the vector $u$ has next coordinates

$$\begin{cases} u_1 = \cos\varphi, \\ u_2 = \sin\varphi, \end{cases} \quad 0 \le \varphi \le \frac{\pi}{2},$$

and the point $A$ lying on the ray produced by vector $u$ has next coordinates

$$\begin{cases} y_{1A} = u_1 \cdot t; \\ y_{2A} = u_2 \cdot t, \end{cases} \quad 0 \le t.$$

Let function $H(y)$ be determined by the formula (7)

$$H(y) = \max_{1 \le i \le k}\{h_i(y)\} \quad (13)$$

# VECTOR OPTIMIZATION BY TWO OBJECTIVE FUNCTIONS

And the angle $\varphi$ is denoted, and then we consider a problem

$$L = t \to \min \qquad (A)$$

where

$$H(u \cdot t) = 0.$$

Let $t_*$ be the solution of the task (A), then it obtains next theorem.

**Theorem.** If the set $Y$ is convex, then this is sufficient for a point $y = u \cdot t_*$ to belong to $Y_*$.

**Proof.** Let us consider the function

$$U(y_1, y_2) = \max\left\{\frac{y_1}{u_1}, \frac{y_2}{u_2}\right\},$$

then the set

$$U(y_1, y_2) \leq t_*$$

is the cone $K$, whose top is at the point $y_{1*} = u_1 \cdot t_*$; $y_{2*} = u_2 \cdot t_*$ (see fig. 6), that lays on the boundary of the set $Y$, and the intersection of this cone with the convex set $Y$ is the point specified. Then by force of the condition (10) we get the proof of the theorem. □

Note, that the task (A) allows to determine the point $y_*$ for non-convex area $Y$, if it satisfies the condition of the area, that we define as a condition (B).

We will explain the condition B by the use of a fig. 6.
1. Any vertical line, that has intersections with the boundary $Y$ has a point, whose second component ($y_2$) is minimal and is not bigger than the second component of a point $A_2$.
2. Any horizontal line that has intersections with the boundary $Y$ has a point whose first component ($y_1$) is minimal and is not bigger than the first component of a point $A_1$.

Or in mathematical terms:

Let $\min y_2$ be the minimal second component of the points of the intersection of vertical line with the boundary of the area $Y$, and $\min y_1$ – minimal first component of the points of the intersection of the horizontal line with the boundary of the area $Y$, then the condition (B) can be formulated as

$$(\min y_1, \min y_2) \notin Y, \quad (\min y_1, \min y_2) \in K_0 \qquad (B)$$

where the cone $K_0$ contains the area $Y$.

**Linear vector optimization problem**

This task for two objective functions is

$$\begin{pmatrix} y_1 \\ y_2 \end{pmatrix} \to \min$$

where

$$Ay \leq b; \quad y \geq 0.$$

As the area $Y$ in this statement of the problem is a convex set, then we can apply the theorem and we get a task like (A).

$$L = t \rightarrow \min$$

$$Au \cdot t - b \leq 0.$$

Thus, when the boundaries are

$$y_1 + y_2 \geq 5;$$
$$y_1 + 3y_2 \geq 8;$$
$$6y_1 + y_2 \geq 14;$$
$$7y_1 + 4y_2 \leq 39,$$

and we add

$$0 \leq y_1 \leq u_1 \cdot t;$$
$$0 \leq y_2 \leq u_2 \cdot t,$$
$$t \geq 0$$

when the vector $u$ we get a simple task of linear programming.

The code of the program for solving this task in the package Maple 7 is

```
>X:=array(1..1000, []);
>Y:= array(1..1000, []);
>k:=0:
>for x0 from 0.01 by 0.01 to 3.14/2 do
    k:=k+1:
    s:={y[1]+y[2]>=5,
        y[1]+3*y[2]>=8,
        6*y[1]+y[2]>=14,
        7*y[1]+4*y[2]<=39,
        y[1]<=cos(x0)*t,
        y[2]<=sin(x0)*t}:
    L:=minimize(L, s, NONNEGATIVE):
    for z in A do
        if op(1,z)<>t then
            if op(1,op(1,z))=1 then X[k]:=op(2,z)
                else Y[k]:=op(2,z):
            end if:
        end if:
    end do:
end do:
>plot([X[j], Y[j], j=1..k], style=point, thickness=3);
```

The result of this program is shown in a fig. 7.

# VECTOR OPTIMIZATION BY TWO OBJECTIVE FUNCTIONS

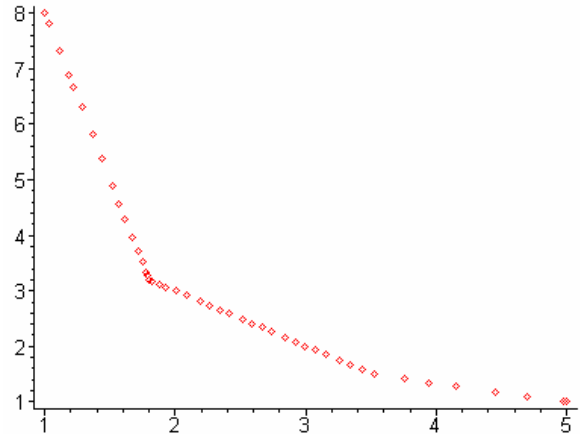

Figure 7. Geometrical interpretation of the solution of the linear vector optimization problem by the program.

As it has pointed the method suggested is close to one of works by A. Messac. Thus, e. g., there is formulated a NBI problem in work [11] that is close to task (A). The main idea of the NBI method is the introduction of the quasi-normal vector $n$. In our task (A) we introduce the vector $u$ that allows to build the cone $K$ and to use the necessary and sufficient condition (10).

In conclusion we show the example from the work [7], where the set $Y$ is not convex, but satisfying the condition (B).

The program of solution of this example in the package Maple 7 is:

```
> X:=array(1..1000,[ ]);Y:=array(1..1000,[ ]);
> h[1]:=1-y[1]^2-y[2]^2/9;
> h[2]:=16-y[1]^4-y[2]^4;
> h[3]:=1-1/27*y[1]^3-y[2]^3;
```

$$h_1 := 1 - y_1^2 - \frac{1}{9} y_2^2$$

$$h_2 := 16 - y_1^4 - y_2^4$$

$$h_3 := 1 - \frac{1}{27} y_1^3 - y_2^3$$

```
> H:=max(h[1], h[2], h[3]);
> k:=0:
> for x0 from 0.1 by 0.01 to 1.47 do
        k:=k+1:
        Hmax:=-10:
        for t from 0.01 by 0.01 to 10 do
                y[1]:=cos(x0)*t:
                y[2]:=sin(x0)*t:
                if H<0 and H>Hmax  then
                        Hmax:=H:
                        tmin:=t:
                end if:
        end do:
        X[k]:=cos(x0)*tmin:
        Y[k]:=sin(x0)*tmin:
end do:
> plot([X[j], Y[j], j=1..k],style=line, thickness=3);
```

The result of the program is shown in a fig. 8

# VECTOR OPTIMIZATION BY TWO OBJECTIVE FUNCTIONS

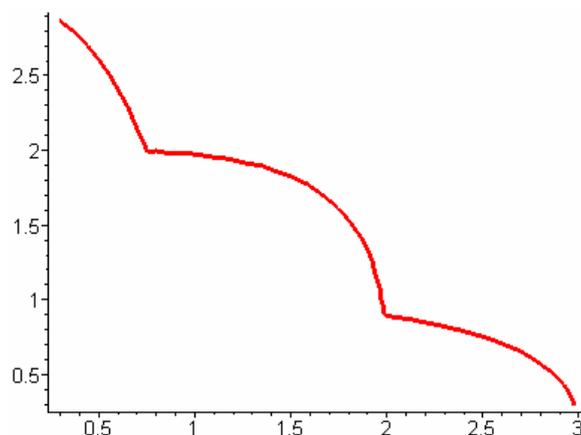

Figure 8. The solution of the vector optimization problem from the work [7] with the use of the task (A)

According to the paper we can conclude:
- If the set $Y$ satisfies the condition (B), then the Pareto`s set $Y_*$ lead to sequence of solving the tasks like the task (A);
- If the set $Y$ does not satisfy the condition (B), then solving the sequence of the tasks (A) we get a set $\tilde{Y}$ that contains $Y_*$.